\documentclass{article}%
\usepackage{amssymb}
\usepackage{amsmath}
\usepackage{amsfonts}
\usepackage{graphicx}%
\providecommand{\U}[1]{\protect\rule{.1in}{.1in}}
\newtheorem{theorem}{Theorem}
\newtheorem{acknowledgement}[theorem]{Acknowledgement}

\newtheorem{corollary}[theorem]{Corollary}

\newtheorem{lemma}[theorem]{Lemma}

\newtheorem{proposition}[theorem]{Proposition}
\newtheorem{remark}[theorem]{Remark}

\begin{document}
The paper will appear in \textbf{Stochastic Processes and Applications}

Available online 17 April 2020

https://doi.org/10.1016/j.spa.2020.04.004

\bigskip

\begin{center}
{\large A new CLT for additive functionals of Markov chains }

\bigskip

Magda Peligrad

\bigskip

Department of Mathematical Sciences, University of Cincinnati, PO Box 210025,
Cincinnati, Oh 45221-0025, USA. \texttt{ }
\end{center}

email: peligrm@ucmail.uc.edu

\bigskip


\smallskip

\bigskip

\begin{center}

Abstract
\end{center}

In this paper we study the central limit theorem for additive functionals of
stationary Markov chains with general state space by using a new idea
involving conditioning with respect to both the past and future of the chain.
Practically, we show that any additive functionals of a stationary and totally
ergodic Markov chain with $\mathrm{var}(S_{n})/n$ uniformly bounded, satisfies
a $\sqrt{n}-$central limit theorem with a random centering. We do not assume
that the Markov chain is irreducible and aperiodic. However, the random
centering is not needed if the Markov chain satisfies stronger forms of
ergodicity. In absence an ergodicity the convergence in distribution still
holds, but the limiting distribution might not be normal.

\section{Introduction}

A basic result in probability theory is the central limit theorem. To go
beyond the independent case, the dependence is often restricted by using
projective criteria. For instance, the martingales are defined by using a
projective condition with respect to the past sigma field. There also is an
abundance of martingale-like conditions, which define classes of processes
satisfying the CLT. Among them Gordin's condition (Gordin, 1969), Gordin and
Lifshits condition (1978), Heyde's projective condition (Heyde, 1974,
Voln\'{y}, 1993), mixingales (McLeish, 1975), Maxwell and Woodroofe condition
(2000), just to name a few. All of them have in common that the conditions are
imposed on the conditional expectation of a variable with respect to the past
sigma field.

There is, however, the following philosophical question. Note that a partial
sum does not depend on the direction of time, i.e.
\[
S_{n}=X_{1}+X_{2}+...+X_{n}=X_{n}+X_{n-1}+...+X_{1}.
\]
However a condition of type "martingale-like" depends on the direction of
time. Therefore, in order to get results for $S_{n},$ it is natural to also
study projective conditions that are symmetric with respect to the direction
of time. Furthermore, many mixing conditions (see Bradley, 2005, for a survey)
and harnesses (see for instance Williams, 1973 and Mansuy and Yor, 2005) are
independent of the direction of time.

For additive functionals of reversible, stationary and ergodic Markov chains,
with centered and square integrable variables, Kipnis and Varadhan (1986)
proved that if $E(S_{n}^{2})/n$ converges to a finite limit, then the CLT
holds. This is not true without assuming reversibility (see for instance
Bradley (1989) or Cuny and Lin (2016), Prop. 9.5(ii), among other examples).
On the other hand, for additive functionals of Harris recurrent and aperiodic
Markov chains with centered and square integrable variables, Chen (1999,
Theorem II. 3.1) proved that if $S_{n}/\sqrt{n}$ is stochastically bounded, it
satisfies the CLT.

These results suggest and motivate the study of limiting distribution for
stationary Markov chains with additive functionals satisfying $\sup_{n}%
E(S_{n}^{2})/n<\infty$. With this aim, we introduce a new idea, which involves
conditioning with respect to both the past and the future of the process. By
using this idea together with a blocking argument and martingale approximation
techniques, we shall prove that functions of a Markov chain which is
stationary and totally ergodic (in the ergodic theoretical sense) satisfy the
CLT,\ provided that we use a random centering and we assume that
$\mathrm{var}(S_{n})/n$ is uniformly bounded. In case when the stationary
Markov chain satisfies stronger forms of ergodicity, the random centering is
not needed. Among these classes are the absolutely regular Markov chains. For
this class, our result gives as a corollary, a new interpretation of the
limiting variance in the CLT in Theorem II. 2.3 of Chen (1999) and a totally
different new approach. We also provide a new proof for the CLT for interlaced
mixing Markov chains. We also point out that, when the Markov chain is
stationary but not necessarily ergodic, the limiting distribution still exists
and we express it as a mixture of distributions.

Our paper is organized as follows. In Section 2 we present the results.
Section 3 is dedicated to their proofs.

\section{Results}

Throughout the paper we assume that $(\xi_{n})_{n\in\mathbb{Z}}$ is a
stationary Markov chain$,$ defined on a probability space $(\Omega
,\mathcal{F},\mathbb{P})$ with values in a measurable space $(S,\mathcal{A})$.
Denote by $\mathcal{F}_{n}=\sigma(\xi_{k},k\leq n)$ and by $\mathcal{F}%
^{n}=\sigma(\xi_{k},k\geq n)$. The marginal distribution on $\mathcal{A}$ is
denoted by $\pi(A)=\mathbb{P}(\xi_{0}\in A)$.{ To settle the concern about the
existence of a Markov chain with general state space,} we shall construct the
Markov chain from a kernel $P(x,A)$, we assume an invariant distribution $\pi$
exists and invoke the Ionescu Tulcea (1949) result.{ }

{Next, let $\mathbb{L}_{0}^{2}(\pi)$ be the set of measurable functions on $S$
such that $\int f^{2}d\pi<\infty$ and $\int fd\pi=0.$ For a function} ${f}\in
${$\mathbb{L}_{0}^{2}(\pi)$ let }%
\begin{equation}
{X_{i}=f(\xi_{i}),\ S_{n}=\sum\nolimits_{i=1}^{n}X_{i}.} \label{def X,S}%
\end{equation}
{ We denote by }${{||X||}}$ the norm in {$\mathbb{L}^{2}$}$(\Omega
,\mathcal{F},\mathbb{P})$.{ }

With the exception of Remark \ref{Rem nonergodic}, in all the other results we
shall assume the total ergodicity of the shift $\theta$ of the sequence
$(\xi_{n})_{n\in\mathbb{Z}}$ with respect to $\mathbb{P},$ i.e. $\theta^{m}$
is ergodic for every $m\geq1$. For the definition of the ergodicity of the
shift we direct the reader to the subsection "A return to Ergodic Theory" in
Billingsley (1995) p. 494.

Let us consider the operator $P$ induced by the kernel $P(x,A)$ on bounded
measurable functions on $(S,\mathcal{A})$ defined by $Pf(x)=\int
\nolimits_{S}f(y)P(x,dy)$. By using Corollary 5 p. 97 in Rosenblatt (1971),
the shift of $(\xi_{n})_{n\in\mathbb{Z}}$ is totally ergodic with respect to
$\mathbb{P}$\ if and only if the powers $P^{m}\ $are ergodic with respect to
$\pi$ for all natural $m$ (i.e. $P^{m}f=f$ for $f$ bounded on $(S,\mathcal{A}%
)$ implies $f$ is constant $\pi-$a.e.). For more information on total
ergodicity, we refer to the survey paper by Quas (2009).

\bigskip

Below, $\Rightarrow$ denotes the convergence in distribution and by
$N(\mu,\sigma^{2})$ we denote a normally distributed random variable with mean
$\mu$ and variance $\sigma^{2}.$

\subsection{Central limit theorem}

We shall establish the following CLT.

\begin{theorem}
\label{Th random center}Assume that
\begin{equation}
\sup_{n\geq1}\frac{E(S_{n}^{2})}{n}<\infty. \label{varsup1}%
\end{equation}
Then, the following limit exists
\begin{equation}
\lim_{n\rightarrow\infty}\frac{1}{n}||S_{n}-E(S_{n}|\xi_{0},\xi_{n}%
)||^{2}=\mathbb{\sigma}^{2} \label{def teta}%
\end{equation}
and%
\[
\frac{S_{n}-E(S_{n}|\xi_{0},\xi_{n})}{\sqrt{n}}\Rightarrow N(0,\mathbb{\sigma
}^{2}).
\]

\end{theorem}

\begin{remark}
It should be noted that, by condition (\ref{varsup1}),\ it follows that
$0\leq\mathbb{\sigma}^{2}<\infty.$ When $\mathbb{\sigma}^{2}=0$ then
$(S_{n}-E(S_{n}|\xi_{0},\xi_{n}))/\sqrt{n}\rightarrow^{\mathbb{P}}0$ as
$n\rightarrow\infty.$ Also, for any stationary sequence, starting from the
identity%
\[
E(S_{n}^{2})=E(X_{0}^{2})+2\sum\nolimits_{k=1}^{n-1}\sum\nolimits_{j=1}%
^{k}E(X_{0}X_{j}),
\]
note that the convergence of sums of the covariances implies that $E(S_{n}%
^{2})/n$ is convergent. Furthermore, if the sums of covariances are bounded by
a constant then (\ref{varsup1}) holds.
\end{remark}

We would like to mention that, as in the stationary martingale case, in the
absence of ergodicity the limiting distribution still exists and it is a
mixture of distributions.

\begin{remark}
\label{Rem nonergodic}If $(\xi_{n})_{n\in\mathbb{Z}}$ is any stationary Markov
chain, $(X_{n})_{n\in\mathbb{Z}}$ and $(S_{n})_{n\geq1}$ are defined by
(\ref{def X,S})\ and (\ref{varsup1}) holds, then there is a random variable
$\eta^{2}$ such that
\begin{equation}
\frac{S_{n}-E(S_{n}|\xi_{0},\xi_{n})}{\sqrt{n}}\Rightarrow\eta^{2}N(0,1),
\label{CLT random}%
\end{equation}
where $\eta^{2}$ is independent of $N(0,1)$.
\end{remark}

\bigskip

If the random centering is not present, we have the following result:

\begin{corollary}
\label{cor CLT}Assume that (\ref{varsup1}) holds and in addition
\[
\frac{E(S_{n}|\xi_{0},\xi_{n})}{\sqrt{n}}\rightarrow^{\mathbb{P}}0\text{ as
}n\rightarrow\infty.
\]
Then the following limit exists
\[
\lim_{n\rightarrow\infty}\frac{\pi}{2n}\left(  E|S_{n}|\right)  ^{2}%
=\mathbb{\sigma}^{2}%
\]
and
\begin{equation}
\frac{S_{n}}{\sqrt{n}}\Rightarrow N(0,\mathbb{\sigma}^{2}). \label{CLT2}%
\end{equation}

\end{corollary}

As an immediate consequence of Theorem \ref{Th random center}, we give next
sufficient conditions for the CLT with the traditional limiting variance.

\begin{corollary}
\label{CorCLTsigma}Assume that (\ref{varsup1}) holds and in addition
\begin{equation}
\lim_{n\rightarrow\infty}\frac{1}{n}||E(S_{n}|\xi_{0},\xi_{n})||^{2}=0.
\label{bad}%
\end{equation}
Then%
\begin{equation}
\lim_{n\rightarrow\infty}\frac{E(S_{n}^{2})}{n}=\sigma^{2} \label{varSn}%
\end{equation}
and
\begin{equation}
\frac{S_{n}}{\sqrt{n}}\Rightarrow N(0,\sigma^{2}). \label{CLT}%
\end{equation}

\end{corollary}

We can also give a sufficient condition for (\ref{bad}) in terms of individual
random variables.

\begin{proposition}
\label{propsufzero}Assume that (\ref{varsup1}) holds and the following
condition is satisfied
\begin{equation}
\lim_{n\rightarrow\infty}n||E(X_{0}|\xi_{-n},\xi_{n})||^{2}=0. \label{badn}%
\end{equation}
Then (\ref{varSn}) and (\ref{CLT}) hold.
\end{proposition}

In Subsection \ref{abs reg}, in the context of absolutely regular Markov
chains, we shall comment that condition (\ref{badn}) alone does not imply
(\ref{CLT}). However, a reinforced condition does:

\begin{corollary}
\label{cor mixi3}Assume that%
\begin{equation}
\sum\nolimits_{k\geq1}||E(X_{0}|\xi_{-k},\xi_{k})||^{2}<\infty.
\label{mixingale}%
\end{equation}
Then (\ref{varSn}) and (\ref{CLT}) hold with $\mathbb{\sigma}^{2}%
=||X_{0}||^{2}+2\sum\nolimits_{k\geq1}E(X_{0}X_{k})$.
\end{corollary}

\subsection{Absolutely regular Markov chains \label{abs reg}}

Relevant to this section is the coefficient of absolute regularity, which was
introduced by Volkonskii and Rozanov (1959) and was attributed there to
Kolmogorov. For a stationary sequence $\mathbf{\xi}=(\xi_{k})_{k\in Z},$ not
necessarily Markov, with values in a separable Banach space, the coefficient
of absolute regularity is defined by
\[
\beta_{n}=\beta_{n}(\mathbf{\xi)=}E\left(  \sup_{A\in\mathcal{F}^{n}%
}|\mathbb{P}(A|\mathcal{F}_{0})-\mathbb{P}(A)|\right)  .
\]
The chain is called absolutely regular if $\beta_{n}\rightarrow0.$ We can
easily see from this definition that $\beta_{n}$ is monotonic. Furthermore
$\beta_{n}$ is symmetric, in the sense that $\beta(\xi_{0},\xi_{n}%
\mathbf{)=}\beta(\xi_{n},\xi_{0}\mathbf{).}$ This fact can be easily seen by
using the equivalent definition for $\beta_{n}$ in terms of partitions (see
Definitions 3.3 and 3.5 in Bradley, 2007). For such sequences, Bradley (1989)
constructed an example of a stationary, pairwise independent, absolutely
regular sequence for which a nondegenrate the central limit theorem cannot hold.

For a Markov chain $\mathbf{\xi}=(\xi_{k})_{k\in Z}$, with values in a
separable Banach space, the coefficient of absolute regularity is equal to
(see Proposition 3.22 (III,5) in Bradley, 2007)
\[
\beta_{n}=\beta_{n}(\mathbf{\xi)=}\beta(\xi_{0},\xi_{n}\mathbf{)=}E\left(
\sup_{A\in\mathcal{B}}|\mathbb{P}(\xi_{n}\in A|\xi_{0})-\mathbb{P}(\xi_{0}\in
A)|\right)  ,
\]
where $\mathcal{B}$ denotes the Borel sigma filed.

Let us mention that there are numerous examples of stationary absolutely
regular Markov chains. For easy reference we refer to Section 3 in Bradley
(2005) survey paper and to the references mentioned there. We know that a
strictly stationary, countable state Markov chain is absolutely regular if and
only if the chain is irreducible and aperiodic. Also, any strictly stationary
Harris recurrent and aperiodic Markov chain is absolutely regular. It is also
well-known that $\beta_{n}\rightarrow0$ implies total ergodicity in the
measure theoretical sense. Also, in many situations these coefficients are
tractable. The computation of the coefficients of absolute regularity is an
area of intense research, with numerous applications to time series and
statistics. There is a vast literature on this subject. See for instance
Davydov (1973), Mokkadem (1990), Doukhan (1994), Doukhan et al. (1994), Ango
Nze (1998), Douc, Moulines and Soulier (2007), Bradley (2007, Vol. 1,2,3)
among others.

Due to their importance for the Monte Carlo simulations, the central limit
theorem for Markov chains was intensively studied under the absolute
regularity condition. In this direction we mention the books by Nummelin
(1984), Meyn and Tweedie (2009) and Chen (1999) and we also refer to the
survey paper by Jones (2004).

In the works mentioned above, the vast majority of results concerning the
CLT\ for absolutely regular Markov chains require sufficient conditions in
terms of moments and mixing rates. Some of them require rates which combine
the tail distribution of a variable with the mixing coefficients.

By using regeneration techniques and partition in independent blocks
(Nummelin's splitting technique, 1978) it was proven that, in this setting, a
necessary and sufficient condition for the CLT is that $S_{n}/\sqrt{n}$ is
stochastically bounded (Theorem II.2.3 in Chen, 1999). However, the limit has
a variance which is described in terms of the split chain and it is difficult
to describe. Our next Corollary is obtained under a more general condition
than in Theorem II. 3.1 in Chen (1999), and sheds new light on the asymptotic
variance in Chen's Theorem II.2.3. The advantage of these results is that no
rate of convergence to zero of the mixing coefficients is required. However,
some information about the variance of partial sums is needed.

\begin{corollary}
\label{Thalpha0}Assume that (\ref{varsup1}) holds and the sequence is
absolutely regular. Then (\ref{CLT2}) holds with $\mathbb{\sigma}^{2}%
=\lim_{n\rightarrow\infty}\pi\left(  E|S_{n}|\right)  ^{2}/2n$.
\end{corollary}

To give a CLT\ where the limiting variance is $\sigma^{2}$ defined in
(\ref{varSn}), we shall verify condition (\ref{badn}) of Proposition
\ref{propsufzero}. Denote by $Q$ the quantile function of $|X_{0}|,$ i.e., the
inverse function of $t\mapsto{\mathbb{P}}(|X_{0}|>t).$ We obtain the following result:

\begin{corollary}
\label{Th alpha}Assume that (\ref{varsup1}) holds and the following condition
is satisfied
\begin{equation}
\lim_{n\rightarrow\infty}n\int\nolimits_{0}^{\beta_{n}}Q^{2}(u)du=0.
\label{cond beta}%
\end{equation}
Then (\ref{varSn}) and (\ref{CLT}) hold.
\end{corollary}

In terms of moments, by H\"{o}lder's inequality, (\ref{cond beta}) is implied
by $E(|X_{0}|^{2+\delta})<\infty$ and $n\beta_{n}^{\delta/(2+\delta
)}\rightarrow0,$ for some $\delta>0$. If $X_{0}$ is bounded a.s., the mixing
rate required for this corollary is $n\beta_{n}\rightarrow0.$

Finally, condition (\ref{mixingale}) is verified if%
\begin{equation}
\sum\nolimits_{n\geq1}\int\nolimits_{0}^{\beta_{n}}Q^{2}(u)du<\infty,
\label{condstrongCLT}%
\end{equation}
and then (\ref{varSn}) and (\ref{CLT}) hold. Further reaching results could be
found in Doukhan et al. (1994), where a larger class of processes was
considered. According to Corollary 1 Doukhan et al. (1994), (\ref{cond beta}%
)\ alone is not enough for (\ref{CLT}). Actually, the stronger condition
(\ref{condstrongCLT}) is a minimal condition for the CLT for $S_{n}/\sqrt{n}$
in the following sense. In their Corollary 1, Doukhan et al. (1994)
constructed a stationary absolutely regular Markov chain $(\xi_{k})_{k\in Z}$
and a function ${f}\in${$\mathbb{L}_{0}^{2}(\pi)$}, which barely does not
satisfies (\ref{condstrongCLT}) and $S_{n}/\sqrt{n}$ does not satisfy the CLT.
For instance, this is the case when for an $a>1,$ $\beta_{n}=cn^{-a}$ and
$Q^{2}(u)$ behaves as $u^{-1+1/a}|\log u|^{-1}$ as $u\rightarrow0^{+}.$ In
this case
\[
\sum\nolimits_{n\geq1}\int\nolimits_{0}^{\beta_{n}}Q^{2}(u)du=\int
\nolimits_{0}^{1}\beta_{n}^{-1}(u)Q^{2}(u)du=\infty
\]
and, according to Corollary 1 in Doukhan et al. (1994), there exists a Markov
chain with these specifications, such that $S_{n}/\sqrt{n}$ does not satisfy
the CLT. However, for these specifications (\ref{cond beta})\ is satisfied. If
in addition we know that (\ref{varsup1}) holds, then, by Corollary
\ref{Th alpha}, the CLT holds for $S_{n}/\sqrt{n}.$

\bigskip

Let us point out for instance, a situation where Corollary \ref{Th alpha} is
useful. Let $\mathbf{Y}=(Y_{i})$ and $\mathbf{Z}=(Z_{j})$ be two absolutely
regular Markov chains of centered, bounded random variables, independent among
them and satisfying the following conditions $\sum\nolimits_{n\geq1}\beta
_{n}(Y)<\infty$ and $n\beta_{n}(Z)\rightarrow0.$ If we define now the sequence
$\mathbf{X}=(X_{n}),$ where for each $n$ we set $X_{n}=Y_{n}Z_{n},$ then, by
Theorem 6.2 in Bradley (2007), we have $\beta_{n}(X)\leq\beta_{n}(Y)+\beta
_{n}(Z).$ As a consequence, by using the monotonicity of $\beta_{n},$ we have
$n\beta_{n}(X)\rightarrow0$ and condition (\ref{cond beta}) is satisfied.
Certainly, this condition alone does not assure that the CLT holds. In order
to apply Corollary \ref{Th alpha} we have to verify that condition
(\ref{varsup1}) holds. Conditioned by $\mathbf{Z}$ the partial sum of
$(X_{n})$ becomes a linear combination of the variables of $\mathbf{Y}$ and we
can apply Corollary 7 in Peligrad and Utev (2006). It follows that there is a
positive constant $C$ such that $\ $%
\begin{equation}
E_{\mathbf{Y}}(\sum\nolimits_{k=1}^{n}Y_{k}Z_{k})^{2}\leq C(\sum
\nolimits_{k=1}^{n}Z_{k}^{2})\text{ a.s.,} \label{varYZ}%
\end{equation}
where $E_{\mathbf{Y}}$ denotes the partial integral with respect to the
variables of $\mathbf{Y}$. By the independence of the sequences $\mathbf{Y}$
and $\mathbf{Z}$ we obtain%
\begin{equation}
E(\sum\nolimits_{k=1}^{n}Y_{k}Z_{k})^{2}\leq CnE(Z_{0}^{2}). \label{CLTWZ}%
\end{equation}
Therefore condition (\ref{varsup1}) holds. By Corollary \ref{Th alpha} we
obtain that (\ref{varSn}) and (\ref{CLT}) hold for the sequence $(X_{n})$.

As a particular example of this kind let us consider two stationary renewal
processes $\mathbf{\xi}=\mathbf{(}\xi_{i}\mathbf{)}$ and $\mathbf{\eta
}=\mathbf{(}\eta_{i}),$ with countable state space $\{0,1,2,...\}$ and
independent among them. For the transition probabilities of $\mathbf{(}\xi
_{i}\mathbf{)}$ we take for $i\geq1,$ $P(\xi_{1}=i-1|\xi_{0}=i)=1,$
$p_{i}=P(\xi_{1}=i|\xi_{0}=0)=\left(  2i^{3}(\log(i+1))^{2}\right)  ^{-1}$,
and $p_{0}=P(\xi_{1}=0|\xi_{0}=0)=1-P(\xi_{1}\geq1|\xi_{0}=0).$

For $\mathbf{(}\eta_{i})$ we take for $i\geq1,$ $P(\eta_{1}=i-1|\eta
_{0}=i)=1,$ $q_{i}=P(\eta_{1}=i|\eta_{0}=0)=\left(  2i^{3}(\log(i+1))\right)
^{-1}$ and $q_{0}=P(\eta_{1}=0|\eta_{0}=0)=1-P(\eta_{1}\geq1|\eta_{0}=0).$
From Theorem 5 in Davydov (1973), we know that the $\beta-$mixing coefficients
for these sequences are of orders%
\[
\beta_{n}(\mathbf{\xi})\leq c\frac{1}{n(\log(n+1))^{2}}\text{ and }\beta
_{n}(\mathbf{\eta})\leq c\frac{1}{n\log(n+1)}\,,
\]
where $c$ is a positive constant. Now, let $f$ and $g$ be two bounded function
and define the sequences $\mathbf{Y}$ and $\mathbf{Z}$ by $Y_{i}=f(\xi
_{i})-E(f(\xi_{i}))$ and $Z_{i}=g(\eta_{i})-E(g(\eta_{i}))$ and set
$X_{i}=X_{i}Y_{i}.\ $Clearly, for this example $\sum\nolimits_{n\geq1}%
\beta_{n}(Y)<\infty$ and $n\beta_{n}(Z)\rightarrow0$ and we can apply
Corollary \ref{Th alpha} for the sequence $(X_{n})$.

\subsection{Interlaced mixing Markov chains}

Another example where Corollary \ref{CorCLTsigma} applies is the class of
interlaced mixing Markov chains. Let $\mathcal{A},\mathcal{B}$ be two sub
$\sigma$-algebras of $\mathcal{F}$. Define the maximal coefficient of
correlation
\[
\rho(\mathcal{A},\mathcal{B})=\sup_{f\in\mathbb{L}_{0}^{2}(\mathcal{A}),\text{
}g\in\mathbb{L}_{0}^{2}(\mathcal{B})}\frac{|E(fg)|}{||f||\cdot||g||}\text{ ,}%
\]
where $\mathbb{L}_{0}^{2}(\mathcal{A})$ ($\mathbb{L}_{0}^{2}(\mathcal{B})$) is
the space of random variables that are $\mathcal{A-}$measurable (respectively
$\mathcal{B-}$measurable$\mathcal{)}$, zero mean and square integrable. For a
sequence of random variables, $(\xi_{k})_{k\in\mathbb{Z}}$, we define
\[
\rho_{n}^{\ast}=\sup\rho(\sigma(\xi_{i},i\in S),\sigma(\xi_{j},j\in T))\text{
,}%
\]
where the supremum is taken over all pairs of disjoint sets, $T$ and $S$ or
$\mathbb{R}$ such that $\min\{|t-s|:t\in T,$ $s\in S\}\geq n.$ We call the
sequence $\rho^{\ast}-$mixing if $\rho_{n}^{\ast}\rightarrow0$ as
$n\rightarrow\infty.$

The $\rho^{\ast}$-mixing condition goes back to Stein (1972) and to Rosenblatt
(1972). It is well-known that $\rho^{\ast}-$mixing implies total ergodicity.
Also, there are known examples (see Example 7.16 in Bradley, 2007) of
$\rho^{\ast}-$mixing sequences which are not absolutely regular.

Our next Corollary shows that our result provides an alternative proof of the
CLT for interlaced $\rho^{\ast}$-mixing Markov chains. Although the result
itself is not new, it provides another example where condition (\ref{bad}) is
verified. For further reaching results see for instance Theorem 11.18 in
Bradley (2007) and Corollary 9.16 in Merlev\`{e}de, Peligrad and Utev (2019).

\begin{corollary}
\label{corinter}Assume that $(\xi_{k})_{k\in\mathbb{Z}}$ is a stationary
$\rho^{\ast}$-mixing Markov chain. Then (\ref{bad}), (\ref{varSn}) and
(\ref{CLT}) hold.
\end{corollary}

\section{Proofs}

\textbf{Proof of Theorem \ref{Th random center}}

\bigskip

The proof of this central limit theorem is based on the martingale
approximation technique. Fix $m$ ($m<n$) a positive integer and make
consecutive blocks of size $m$. Denote by $Y_{k}$ the sum of variables in the
$k$'th block. Let $u=u_{n}(m)=[n/m].$ So, for $k=0,1,...,u-1,$ we have%
\[
Y_{k}=Y_{k}(m)=(X_{km+1}+...+X_{(k+1)m}).
\]
Also denote%
\[
Y_{u}=Y_{u}(m)=(X_{um+1}+...+X_{n}).
\]
For $k=0,1,...,u-1$ let us consider the random variables\textbf{ }%
\[
D_{k}=D_{k}(m)=\frac{1}{\sqrt{m}}(Y_{k}-E(Y_{k}|\xi_{km},\xi_{(k+1)m})).
\]
By the Markov property, conditioning by $\sigma(\xi_{km},\xi_{(k+1)m})$ is
equivalent to conditioning by $\mathcal{F}_{km}\vee\mathcal{F}^{(k+1)m}.$ Note
that $D_{k}$ is adapted to $\mathcal{F}_{(k+1)m}=\mathcal{G}_{k}.$ Then we
have $E(D_{1}|\mathcal{G}_{0})=0$ a.s. Since we assumed that the shift
$\theta$ of the sequence $(\xi_{n})_{n\in\mathbb{Z}}$ is totally ergodic, we
deduce that we have a stationary and ergodic sequence of square integrable
martingale differences $(D_{k},\mathcal{G}_{k})_{k\geq0}$.

Therefore, by the classical central limit theorem for ergodic martingales, for
every $m,$ a fixed positive integer, we have
\[
\frac{1}{\sqrt{u}}M_{u}(m):=\frac{1}{\sqrt{u}}\sum\nolimits_{k=0}^{u-1}%
D_{k}(m)\Rightarrow N_{m}\text{ as }n\rightarrow\infty,
\]
where $N_{m}$ is a normally distributed random variable with mean $0$ and
variance $m^{-1}||Y_{0}-E(Y_{0}|\xi_{0},\xi_{m})||^{2}$. Now consider
$(m^{\prime})$ a subsequence of $\mathbb{N}$ such that
\begin{equation}
\lim_{m^{\prime}\rightarrow\infty}\frac{1}{m^{\prime}}||Y_{0}-E(Y_{0}|\xi
_{0},\xi_{m^{\prime}})||^{2}=\lim\sup_{n\rightarrow\infty}\frac{1}{n}%
||S_{n}-E(S_{n}|\xi_{0},\xi_{n})||^{2}=\eta^{2}. \label{limsup eta}%
\end{equation}
Note that by (\ref{varsup1})\ it follows that $\eta^{2}<\infty.$ This means
that%
\[
N_{m^{\prime}}\Rightarrow N(0,\eta^{2})\text{ as }m^{\prime}\rightarrow
\infty.
\]
Whence, according to Theorem 3.2 in Billingsley (1999),\ in order to establish
the CLT for $\sum\nolimits_{i=0}^{u-1}Y_{i}/\sqrt{n}$ we have only to show
that
\begin{equation}
\lim_{m^{\prime}\rightarrow\infty}\lim\sup_{n\rightarrow\infty}||\frac
{1}{\sqrt{n}}\left(  S_{n}-E(S_{n}|\xi_{0},\xi_{n})\right)  -\frac{1}{\sqrt
{u}}M_{u}(m^{\prime})||^{2}=0. \label{negli1}%
\end{equation}
Denote by $Z_{k}=m^{-1/2}E(Y_{k}|\xi_{km},\xi_{(k+1)m})$ and let
$R_{u}(m)=\sum\nolimits_{k=0}^{u-1}Z_{k}.$ Set
\begin{equation}
S_{u}(m)=M_{u}(m)+R_{u}(m). \label{decom}%
\end{equation}
Let us show that $M_{n}(m)$ and $R_{n}(m)$ are orthogonal. We show this by
analyzing the expected value of all the terms of the product $M_{n}%
(m)R_{n}(m)$. Note that if $j<k,$ since $\mathcal{F}_{(j+1)m}\subset
\mathcal{F}_{km},$ we have
\begin{align*}
&  E[(Y_{k}-E(Y_{k}|\xi_{km},\xi_{(k+1)m}))E(Y_{j}|\xi_{jm},\xi_{(j+1)m})]\\
&  =E[E(Y_{k}-E(Y_{k}|\mathcal{F}_{km}\vee\mathcal{F}^{(k+1)m})|\mathcal{F}%
_{(j+1)m})E(Y_{j}|\xi_{jm},\xi_{(j+1)m})]=0.
\end{align*}
On the other hand, if $j>k,$ since $\mathcal{F}^{_{jm}}\subset\mathcal{F}%
^{(k+1)m}$ then
\begin{align*}
&  E[(Y_{k}-E(Y_{k}|\xi_{km},\xi_{(k+1)m}))E(Y_{j}|\xi_{jm},\xi_{(j+1)m})]\\
&  =E[E(Y_{k}-E(Y_{k}|\mathcal{F}_{km}\vee\mathcal{F}^{(k+1)m})|\mathcal{F}%
^{_{jm}})E(Y_{j}|\xi_{jm},\xi_{(j+1)m})]=0.
\end{align*}
For $j=k,$ by conditioning with respect to $\sigma(\xi_{km},\xi_{(k+1)m}),$ we
note that
\[
E[(Y_{k}-E(Y_{k}|\xi_{km},\xi_{(k+1)m}))E(Y_{k}|\xi_{km},\xi_{(k+1)m})]=0.
\]
Therefore $M_{n}(m)$ and $R_{n}(m)$ are indeed orthogonal. By using now the
decomposition (\ref{decom}), the fact that $M_{n}(m)$ and $R_{n}(m)$ are
orthogonal and $M_{n}(m)$ is a martingale, we obtain the identity
\begin{equation}
\frac{1}{u}||S_{u}(m)||^{2}=\frac{1}{m}||S_{m}-E(S_{m}|\xi_{0},\xi_{m}%
)||^{2}+\frac{1}{u}||R_{u}(m)||^{2}. \label{ID}%
\end{equation}
Also, note that (\ref{varsup1}) and the definition of $Y_{u}$ imply that for
some positive constant $C,$ we have $||Y_{u}||\leq Cm$. Hence, by the
properties of conditional expectations, for every $m$ fixed, we have
\begin{equation}
||\frac{1}{\sqrt{n}}\left(  S_{n}-E(S_{n}|\xi_{0},\xi_{n})\right)  -\frac
{1}{\sqrt{u}}(S_{u}(m)-E(S_{u}(m)|\xi_{0},\xi_{n}))||\rightarrow0\text{ as
}n\rightarrow\infty. \label{neglast}%
\end{equation}
Recall the definition of $M_{u}(m)$, which is orthogonal to $\mathcal{F}%
_{0}\vee\mathcal{F}^{um}.$ By using again the properties of conditional
expectations and the identity (\ref{ID}), for every $m$ we have%
\begin{gather*}
\frac{1}{u}||S_{u}(m)-E(S_{u}(m)|\xi_{0},\xi_{n})-M_{u}(m)||^{2}=\frac{1}%
{u}||R_{u}(m)||^{2}-\frac{1}{u}||E(S_{u}(m)|\xi_{0},\xi_{n})||^{2}\\
=\frac{1}{u}\left(  ||S_{u}(m)||^{2}-||E(S_{u}(m)|\xi_{0},\xi_{n}%
)||^{2}\right)  -\frac{1}{m}||S_{m}-E(S_{m}|\xi_{0},\xi_{m})||^{2}\\
=\frac{1}{u}||(S_{u}(m)-E(S_{u}(m)|\xi_{0},\xi_{n})||^{2}-\frac{1}{m}%
||S_{m}-E(S_{m}|\xi_{0},\xi_{m})||^{2}.
\end{gather*}
By passing now to the limit in the last identity with $n\rightarrow\infty,$ by
(\ref{limsup eta}) and (\ref{neglast}) we obtain
\begin{align*}
&  \lim\sup_{n\rightarrow\infty}||\frac{1}{\sqrt{n}}\left(  S_{n}-E(S_{n}%
|\xi_{0},\xi_{n})\right)  -\frac{1}{\sqrt{u}}M_{u}(m)||^{2}\\
&  =\eta^{2}-\left(  \frac{1}{m}||S_{m}-E(S_{m}|\xi_{0},\xi_{m})||^{2}\right)
.
\end{align*}
By letting now $m^{\prime}\rightarrow\infty$ on the subsequence defined in
(\ref{limsup eta}) and taking into account (\ref{neglast}), we have that
(\ref{negli1}) follows. Therefore%
\[
\frac{S_{n}-E(S_{n}|\xi_{0},\xi_{n})}{\sqrt{n}}\Rightarrow N(0,\eta^{2}).
\]
Then, by (\ref{limsup eta}), Skorohod's representation theorem (i.e. Theorem
6.7 in Billingsley, 1999) and by Fatou's lemma we get%
\[
\lim\sup_{n\rightarrow\infty}\frac{E(S_{n}-E(S_{n}|\xi_{0},\xi_{n}))^{2}}%
{n}=\eta^{2}\leq\lim\inf_{n\rightarrow\infty}\frac{E(S_{n}-E(S_{n}|\xi_{0}%
,\xi_{n}))^{2}}{n}.
\]
It follows that (\ref{def teta}) holds as well as the CLT in Theorem
\ref{Th random center}. $\ \square$

\bigskip

\textbf{Proof of Corollary \ref{cor CLT}}.

\bigskip

From Theorem \ref{Th random center} and Theorem 3.1 in Billingsley (1999) we
immediately obtain (\ref{CLT2}). Note that, by (\ref{varsup1}), we have that
$|S_{n}|/\sqrt{n}$ is uniformly integrable and therefore, by (\ref{CLT2}) and
the convergence of moments theorem (Theorem 3.5, Billingsley, 1999)\ we have
that $E|S_{n}|/\sqrt{n}\rightarrow\sqrt{2/\pi}\mathbb{\sigma}$. \ $\square$

\bigskip

\textbf{Proof of Remark \ref{Rem nonergodic}}

\bigskip

The proof of this remark is based on two facts.

\textbf{Fact 1}. Raikov-type CLT for stationary martingale differences. (see
Theorem 3.6 in Hall and Heyde, 1980). If $(D_{k})_{k\in\mathbb{Z}}$ is a
square integrable sequence of martingale differences and $M_{n}=D_{1}%
+...+D_{n},$ then there is a random variable $\eta^{2}$ such that
\[
\frac{M_{n}}{\sqrt{n}}\Rightarrow\eta^{2}N(0,1),
\]
where $\eta^{2}$ is independent on $N(0,1).$

\textbf{Fact 2}. A variant of Theorem 3.2 in Billingsley (1999) for complete
separable metric spaces (Theorem 2 in Dehling et al., 2009). For random
variables $(X_{n}(m^{\prime}),Y_{n})$ with $n\in\mathbb{N}$ and $m^{\prime}$
belonging to a subsequence of $\mathbb{N}$ which tends to $\infty,$ assume
that for every $\varepsilon>0$%
\[
\lim_{m^{\prime}\rightarrow\infty}\lim\sup_{n\rightarrow\infty}P(|X_{n}%
(m^{\prime})-Y_{n}|>\varepsilon)=0,
\]
and for every $m^{\prime},$ $X_{n}(m^{\prime})\Rightarrow Z(m^{\prime})$ as
$n\rightarrow\infty.$ Then there is a random variable $X$ such that
$Z(m^{\prime})\Rightarrow X$ as $m^{\prime}\rightarrow\infty$ and
$Y_{n}\Rightarrow X$ as $n\rightarrow\infty.$

To prove Remark \textbf{\ref{Rem nonergodic}, }we define the subsequence
$(m^{\prime})$ by (\ref{limsup eta})\ and start from relation (\ref{negli1}).
We apply next Fact 1 to the sequence of stationary martingale differences
$(D_{k}(m^{\prime}))_{k\geq0}$ and obtain that
\[
\frac{M_{u}(m^{\prime})}{\sqrt{m^{\prime}}}\Rightarrow\eta_{m^{\prime}}%
^{2}N(0,1)\text{ as }m^{\prime}\rightarrow\infty,
\]
where $\eta_{m^{\prime}}^{2}$ are random variables independent on $N(0,1).$
Then, we apply Fact 2 and deduce that, for some random variable $X,$ both
$\eta_{m^{\prime}}N(0,1)\Rightarrow X$ and $S_{n}/\sqrt{n}\Rightarrow X.$ But
the characteristic function of $\eta_{m^{\prime}}^{2}N(0,1)$ is $E\left(
\exp(-t^{2}\eta_{m^{\prime}}^{2}/2)\right)  $ and therefore $\eta_{m^{\prime}%
}^{2}$ is converges in distribution to some random variable $\eta$ implying
(\ref{CLT random}). \ $\square$

\textbf{\bigskip}

\textbf{Proof of Proposition \ref{propsufzero}}

\bigskip

This proposition follows by applying Corollary \ref{CorCLTsigma}. Note that we
have only to show that (\ref{badn}) implies (\ref{bad}).

We start the proof of this fact by fixing $0<\varepsilon<1$ and writing
\[
S_{n}=S_{[\varepsilon n]}+V_{n}(\varepsilon)+(S_{n}-S_{n-[\varepsilon n]}),
\]
where%
\[
V_{n}(\varepsilon)=\sum\nolimits_{j=[\varepsilon n]+1}^{n-[\varepsilon
n]}X_{j}.
\]
Note that, by the triangle inequality, properties of the norm of the
conditional expectation, condition (\ref{varsup1}) and stationarity, we easily
get
\begin{equation}
\lim\sup_{n\rightarrow\infty}\frac{1}{\sqrt{n}}||E(S_{n}|\xi_{0},\xi
_{n})||\leq\lim_{\varepsilon\rightarrow0}\lim\sup_{n\rightarrow\infty}\frac
{1}{\sqrt{n}}||E(V_{n}(\varepsilon)|\xi_{0},\xi_{n})||. \label{ineqlim}%
\end{equation}
By the Cauchy-Schwartz inequality, for $1\leq a\leq b\leq n,$%
\[
||E(\sum\nolimits_{j=a}^{b}X_{j}|\xi_{0},\xi_{n})||^{2}\leq n\sum
\nolimits_{j=a}^{b}||E(X_{j}|\xi_{0},\xi_{n})||^{2}.
\]
So, by stationarity
\[
\frac{1}{n}||E(V_{n}(\varepsilon)|\xi_{0},\xi_{n})||^{2}\leq\sum
\nolimits_{j=[\varepsilon n]+1}^{n-[\varepsilon n]}||E(X_{0}|\xi_{-j}%
,\xi_{n-j})||^{2}.
\]
Since for $[\varepsilon n]+1\leq j\leq n-[\varepsilon n]$ we have
$\mathcal{F}_{-j}\vee\mathcal{F}^{n-j}\subset\mathcal{F}_{-[\varepsilon
n]}\vee\mathcal{F}^{[\varepsilon n]}$ it follows that%
\[
\frac{1}{n}||E(V_{n}(\varepsilon)|\xi_{0},\xi_{n})||^{2}\leq n||E(X_{0}%
|\xi_{-[\varepsilon n]},\xi_{\lbrack\varepsilon n]})||^{2}.
\]
We obtain (\ref{bad}) by passing to the limit with $n\rightarrow\infty$ in the
last inequality and taking into account (\ref{badn}) and (\ref{ineqlim}).
\ $\square$

\bigskip

\textbf{Proof of Corollary \ref{cor mixi3}}

\bigskip

By the monotonicity of $||E(X_{0}|\xi_{-k},\xi_{k})||$, condition
(\ref{mixingale}) implies (\ref{badn}). Condition (\ref{mixingale}) also
implies the couple of conditions%
\begin{equation}
\sum\nolimits_{k\geq1}||E(X_{0}|\xi_{-k})||^{2}<\infty\text{ and }%
\sum\nolimits_{k\geq1}||E(X_{0}|\xi_{k})||^{2}<\infty. \label{forvar}%
\end{equation}
Now note that by the properties of the conditional expectations, the Markov
property and stationarity, for all $k\geq1$ we easily obtain%
\begin{gather*}
|E(X_{0}X_{2k})|=|E(X_{0}E(X_{2k}|\xi_{k}))|=|E(E(X_{0}|\xi_{k})E(X_{2k}%
|\xi_{k}))|\leq\\
||E(X_{0}|\xi_{k})||\cdot||E(X_{0}|\xi_{-k})||\leq(||E(X_{0}|\xi_{-k}%
)||^{2}+||E(X_{0}|\xi_{k})||^{2})/2.
\end{gather*}
A \ similar relation holds for $|E(X_{0}X_{2k+1})|.$ Hence the two conditions
in (\ref{forvar})\ lead to (\ref{varsup1}). The result follows by applying
Proposition \textbf{\ref{propsufzero}.}\ $\square$

\bigskip

Before proving the corollaries in Subsection \ref{abs reg} we give a more
general definition of the coefficient of absolute regularity $\beta.$ As in
relation (5) in Proposition 3.22 in Bradley, given two sigma algebras
$\mathcal{A}$ and $\mathcal{B}$ with $\mathcal{B}$ separable and for any
$B\in\mathcal{B}$ there is a regular conditional probability $P(B|\mathcal{A}%
)$, then%
\[
\beta(\mathcal{A},\mathcal{B})=E(\sup_{B\in\mathcal{B}}|P(B|\mathcal{A}%
)-P(B)|).
\]
We also need a technical lemma whose proof is given later.

\begin{lemma}
\label{strong}Let $X,Z$ be two random variables on a probability space
$(\Omega,\mathcal{K},P)$ with values in a separable Banach space. Let
$\mathcal{B\subset K}$ be a sub $\sigma-$algebra. Assume that $X$ and $Z$ are
conditionally independent given $\mathcal{B}$. Then
\[
\beta(\mathcal{B},\mathcal{A}\vee\mathcal{C)\leq}\beta(\mathcal{A}%
,\mathcal{B)+}\beta(\mathcal{C},\mathcal{B)+}\beta(\mathcal{A},\mathcal{C)}%
\text{,}%
\]
where $\mathcal{A=\sigma(}X)$ and $\mathcal{C=\sigma(}Z).$
\end{lemma}

\bigskip

\textbf{Proof of Corollary \ref{Thalpha0}}

\bigskip

In order to apply Corollary \ref{cor CLT}\ it is enough to verify that
\[
\frac{E|E(S_{n}|\xi_{0},\xi_{n})|}{\sqrt{n}}\rightarrow0.
\]
Let $v\leq n$ be a positive integer. Then, by (\ref{varsup1}) we have
\begin{equation}
\lim\sup_{n\rightarrow\infty}\frac{E|E(S_{n}|\xi_{0},\xi_{n})|}{\sqrt{n}}%
=\lim\sup_{n\rightarrow\infty}\frac{E|E(V_{n}(v)|\xi_{0},\xi_{n})|}{\sqrt{n}},
\label{ineq V}%
\end{equation}
where $V_{n}(v)=\sum\nolimits_{j=v+1}^{n-v}X_{j}.$ But, it is well-known that
(see Ch.4 in Bradley 2007)%
\[
E|E(V_{n}(v)|\xi_{0},\xi_{n})|\leq8\beta^{1/2}(\sigma(\xi_{i};v\leq i\leq
n-v),\sigma(\xi_{0},\xi_{n}))||S_{n-2v}||_{2}.
\]
By Lemma \ref{strong}, applied with $\mathcal{A=}\sigma(\xi_{0})$,
$\mathcal{B=}\sigma(\xi_{i};$ $v\leq i\leq n-v)$, and $\mathcal{C=}\sigma
(\xi_{n})$ and taking into account the properties $\beta_{v}$ listed at the
beginning of Subsection \ref{abs reg} along with stationarity, we obtain that
\[
\beta(\sigma(\xi_{i};m\leq i\leq n-m),\sigma(\xi_{0},\xi_{n}))\leq\beta
(\xi_{0},\xi_{v})+\beta(\xi_{n},\xi_{n-v})+\beta(\xi_{0},\xi_{n})\leq
3\beta_{v}.
\]
Therefore, for all $v\in\mathbb{N}$
\[
\lim\sup_{n\rightarrow\infty}\frac{E|E(S_{n}|\xi_{0},\xi_{n})|}{\sqrt{n}}%
\leq24\beta_{v}^{1/2}\sup_{n}\frac{1}{\sqrt{n}}||S_{n}||_{2},
\]
and the result follows by letting $v\rightarrow\infty.$ \ $\square$

\bigskip

\textbf{Proof of Corollary \ref{Th alpha}}

\bigskip

This corollary follows by verifying the conditions of Proposition
\ref{propsufzero}. By Rio's (1993) covariance inequality (see also Theorem 1.1
in Rio 2017) we know that%
\[
||E(X_{0}|\xi_{-n},\xi_{n})||^{2}\leq2\int\nolimits_{0}^{\bar{\beta}_{n}}%
Q^{2}(u)du,
\]
where $\bar{\beta}_{n}=\beta(\sigma(\xi_{0}),\sigma(\xi_{-n},\xi_{n})).$

But, according to Lemma \ref{strong}, applied with $\mathcal{A=}\sigma
(\xi_{-n})$, $\mathcal{B=}\sigma(\xi_{0})$, and $\mathcal{C=}\sigma(\xi_{n})$
we obtain
\[
\bar{\beta}_{n}=\beta(\sigma(\xi_{0}),\sigma(\xi_{-n},\xi_{n}))\leq
3\beta(\sigma(\xi_{0}),\sigma(\xi_{n}))=3\beta_{n}%
\]
\ and the result follows. $\ \square$

\bigskip

\textbf{Proof of Corollary \ref{corinter}}

\bigskip

For this class of random variables it is well-known that condition
(\ref{varsup1}) is satisfied (see for instance Lemma 8.23 in Bradley, 2007).
According to Corollary (\ref{CorCLTsigma})\ we have only to verify condition
(\ref{bad}). Note that by (\ref{ineq V})\ it is enough to show that%
\[
\lim_{v\rightarrow\infty}\lim\sup_{n\rightarrow\infty}\frac{1}{\sqrt{n}%
}||E(V_{n}(v)|\xi_{0},\xi_{n})||=0,
\]
with $V_{n}(v)=\sum\nolimits_{j=v+1}^{n-v}X_{j}.$ By the definition of
$\rho_{v}^{\ast}$ we observe that
\[
||E(V_{n}(v)|\xi_{0},\xi_{n})||^{2}=|E(V_{n}(v)E(V_{n}(\varepsilon)|\xi
_{0},\xi_{n})|\leq\rho_{v}^{\ast}||V_{n}(v)||\cdot||E(V_{n}(v)|\xi_{0},\xi
_{n})||.
\]
Whence,
\[
\frac{1}{\sqrt{n}}||E(V_{n}(v)|\xi_{0},\xi_{n})||\leq\rho_{v}^{\ast}\frac
{1}{\sqrt{n}}||V_{n}(v)||.
\]
The result follows by (\ref{varsup1}). $\square$

\bigskip

\bigskip

\textbf{Proof of Lemma \ref{strong}}

\bigskip

Denote the law of $X$ by $P_{X},$ the law of $Z$ by $P_{Z}.$ Also by
$P_{X|\mathcal{B}}$ we denote the regular conditional distribution of $X$
given $\mathcal{B}$ and by $P_{Z|\mathcal{B}}$ the regular conditional
distribution of $Z$ given $\mathcal{B}$. By using the definition of
$\beta(\mathcal{B},\mathcal{A}\vee\mathcal{C)}$ we have to evaluate the
expression $I=E[\sup_{H}|P(H|\mathcal{B})-P(H)|],$ where the supremum is taken
over all $H\subset\mathcal{A}\vee\mathcal{C}$. Denote by $I_{H}$ the indicator
function of $H.$ Since $X$ and $Z$ are conditionally independent given
$\mathcal{B}$ we have
\[
P(H|\mathcal{B})=\iint I_{H}(x,z)P_{X|\mathcal{B}}(dx)P_{Z|\mathcal{B}%
}(dz)\text{ a.s}.
\]
Also,
\[
P(H)=\iint I_{H}(x,z)P_{(X,Z)}(dx,dz).
\]
By the triangle inequality we can write $I\leq I_{1}+I_{2}+I_{3}$ where
\[
I_{1}=E\left(  \sup_{H}|\iint I_{H}(x,z)\left(  P_{X|\mathcal{B}%
}(dx)P_{Z|\mathcal{B}}(dz)-P_{X}(dx)P_{Z|\mathcal{B}}(dz)\right)  |\right)  ,
\]%
\[
I_{2}=E\left(  \sup_{H}|\iint I_{H}(x,z)\left(  P_{X}(dx)P_{Z|\mathcal{B}%
}(dz)-P_{X}(dx)P_{Z}(dz)\right)  |\right)
\]
and%
\[
I_{3}=\sup_{H}|\iint I_{H}(x,z)\left(  P_{X,Z}(dx,dz)-P_{X}(dx)P_{Z}%
(dz)\right)  |.
\]
Now, because$\ I_{H}$ is bounded by $1$, we get%
\[
I_{1}\leq E\left(  \sup_{D\subset\mathcal{R}}|P_{X|\mathcal{B}}(D)-P_{X}%
(D)|\right)  ,
\]%
\[
I_{2}\leq E\left(  \sup_{D\subset\mathcal{R}}|P_{Z|\mathcal{B}}(D)-P_{Z}%
(D)|\right)  ,
\]
and%
\[
I_{3}\leq\iint|P_{X,Z}(dx,dz)-P_{X,Z^{\ast}}|dxdz,
\]
where $\mathcal{R}$ denote the Borel sigma field and $Z^{\ast}$ is a random
variable distributed as $Z$ and independent of $X.$

The result follows by using the definition of $\beta$ and Theorem 3.29, both
in Bradley (2007). $\square$

\bigskip

\begin{acknowledgement}
This paper was partially supported by the NSF grant DMS-1811373. The author
would like to thank Christophe Cuny for discussions about the variance of
partial sums, to Wlodek Bryc for discussions about double tail sigma field and
harnesses and to Richard Bradley for correcting and suggesting a proof of
Lemma \ref{strong}. Special thanks go to the two anonymous referees for
carefully reading of the manuscript and for numerous suggestions, that
improved the presentation of this paper.
\end{acknowledgement}


\begin{thebibliography}{99}                                                                                               %


\bibitem {A}Ango Nze, P. (1998). Crit\`{e}res d'ergodicit\'{e}
g\'{e}om\'{e}trique ou arithm\'{e}tique de mod\`{e}les lin\'{e}aires
perturb\'{e}s \`{a} repr\'{e}sentation markovienne. Comptes Rendus de
l'Acad\'{e}mie des Sciences. Paris 326, S\'{e}rie 1, 371--376.

\bibitem {B95}Billingsley, P. (1995). Probability and measure. Third edition,
Wiley, New York.

\bibitem {B}Billingsley, P. (1999). Convergence of probability measures.
Second edition. Wiley, New York.

\bibitem {Br89}Bradley, R.C. (1989). Stationary, pairwise independent,
absolutely regular sequence for which the central limit theorem fails. Probab.
Theory and Related Fields 81 1-10.

\bibitem {BR0}Bradley, R.C. (2005). Basic properties of strong mixing
conditions. A survey and some open questions. Probability Surveys 2 107-144.

\bibitem {BR}Bradley, R.C. (2007).\textit{ }Introduction to strong mixing
conditions 1, 2, 3. Kendrick Press, Heber City, UT.

\bibitem {Chen}Chen, X. (1999). Limit theorems for functionals of ergodic
Markov chains with general state space. Memoirs of the American Mathematical
Society 139.

\bibitem {CL}Cuny, C. and Lin, M. (2016). Limit theorems for Markov chains by
the symmetrization method. J. Math. Anal. Appl. 434 52--83.

\bibitem {Da}Dehling, H., Durieu, O. and Voln\'{y}, D. (2009). New techniques
for empirical processes of dependent data. Stochastic Process. Appl. 119 3699--3718.

\bibitem {Dav}Davydov, Yu. A. (1973). Mixing conditions for Markov chains.
Theory Probab. Appl. 18 312-328.

\bibitem {DOUC}Douc, R., Moulines, E. and Soulier, P. (2007). Computable
convergence rates for subgeometric ergodic Markov chains. Bernoulli 13 831--848.

\bibitem {DOU}Doukhan, P. (1994). Mixing: Properties and Examples. Lecture
Notes in Statistics 85. Springer Verlag.

\bibitem {Douk}Doukhan, P., Massart, P. and Rio, E. (1994). The functional
central limit theorem for strongly mixing processes. Ann. Inst. H.
Poincar\'{e} Probab. Statist. 30 63-82.

\bibitem {G69}Gordin, M. I. (1969). The central limit theorem for stationary
processes, Soviet. Math. Dokl. 10 1174--1176.

\bibitem {GL}Gordin, M. and Lifshits, B. (1978). The central limit theorem for
stationary Markov processes. Soviet. Math. Dokl. 19 392-393.

\bibitem {HH}Hall, P. and Heyde, C.C. (1980). Martingale limit theory and its
application. Academic Press, New York.

\bibitem {He}Heyde, C. C. (1974). On the central limit theorem for stationary
processes. Z. Wahrsch. Verw. Gebiete 30 315-320.

\bibitem {IT}Ionescu Tulcea, C. T. (1949). Mesures dans les espaces produits.
Atti Accad. Naz. Lincei Rend. 7 208--211.

\bibitem {Jo}Jones, G. L. (2004). On the Markov chain central limit theorem.
Probab. Surv. 1 299--320.

\bibitem {KV}Kipnis, C. and Varadhan, S.R.S. (1986). Central limit theorem for
additive functionals of reversible Markov processes and applications to simple
exclusions. Comm. Math. Phys. 104 1-19.

\bibitem {MY}Mansuy, R. and Yor, M. (2005). Harnesses, L\'{e}vy bridges and
Monsieur Jourdain. Stochastic Process. Appl. 115 329--338.

\bibitem {MPU}Merlev\`{e}de, F., Peligrad, M. and Utev, S. (2019). Functional
Gaussian approximation for dependent structures. Oxford University Press.

\bibitem {Mok}Mokkadem, A. (1990). Propri\'{e}t\'{e}s de m\'{e}lange des
mod\'{e}les autor\'{e}gressifs polynomiaux. Ann. I. H. P. ser. B 26 (2) 219-260.

\bibitem {MW}Maxwell, M. and Woodroofe, M. (2000). Central limit theorems for
additive functionals of Markov chains. Ann. Probab. 28 713--724.

\bibitem {Mc}McLeish, D. L. (1975). A generalization of martingales and mixing
sequences. Adv. in Appl. Probab. 7 247-258.

\bibitem {MT}Meyn S.P. and Tweedie R.L. (2009). Markov chains and stochastic
stability, 2nd edition. Cambridge University Press.

\bibitem {N}Nummelin, E. (1978). A splitting technique for Harris recurrent
chains. Z. Wahrs. verw Gebiete 43 309-318.

\bibitem {Nb}Nummelin, E. (1984). General irreducible Markov chains and
non-negative operators. Cambridge University Press, Cambridge, England.

\bibitem {PU}Peligrad, M. and Utev, S. (2006). Central limit theorem for
stationary linear processes. Ann. Probab. 34 1608--1622.

\bibitem {Q}Quas, A. (2009). Ergodicity and mixing properties, in: R. E.
Meyers (ed.), Encyclopedia of Complexity and Systems Science 2918--2933, Springer.

\bibitem {Rio93}Rio, E. (1993). Covariance inequalities for strongly mixing
processes. Ann. Inst. Henri Poincar\'{e} Probab. Stat. 29 587--597.

\bibitem {Rio17}Rio, E. (2017). Asymptotic theory of weakly dependent random
processes. Springer.

\bibitem {Ros}Rosenblatt, M. (1971). Markov processes. Structure and
asymptotic behavior. Springer, Berlin.

\bibitem {Ros2}Rosenblatt, M. (1972). Central limit theorems for stationary
processes. Proceedings of the Sixth Berkeley Symposium on Probability and
Statistics, vol. 2, 551-561. University of California Press, Los Angeles.

\bibitem {Stein}Stein C. (1972). A bound for the error in the normal
approximation to the distribution of a sum of dependent random variables. In:
Proceedings of the Sixth Berkeley Symposium on Probability and Statistics,
Vol. 2, 583-602. University of California Press, Los Angeles.

\bibitem {VR}Volkonskii, V.A. and Rozanov Yu.A. (1959). Some limit theorems
for random functions I. Theor. Probab. Appl. 4 178-197.

\bibitem {Vol}Voln\'{y}, D. (1993). Approximating martingales and the central
limit theorem for strictly stationary processes. Stochastic Processes and
their Applications 44 (1993) 41-74.

\bibitem {wi}Williams, D. (1973). Some basic theorems on harnesses. Stochastic
Analysis (a tribute to the memory of Rollo Davidson), 349--363. Wiley, London.\ 
\end{thebibliography}
\end{document}